\pdfoutput=1
\documentclass{article}
\usepackage[utf8]{inputenc}
\usepackage{graphicx}
\usepackage{url,xcolor,xspace}
\usepackage{amsmath}
\usepackage{comment}
\usepackage{hyperref}
\usepackage[export]{adjustbox}
\usepackage{geometry}[margin=0.25in]
\usepackage{enumerate}
\usepackage[shortlabels]{enumitem}
\usepackage{amsfonts}
\usepackage{multicol}
\title{Control of Line Pack in Natural Gas System: Balancing Limited Resources under Uncertainty}
\title{\Large  \bf Control of Line Pack in Natural Gas System: Balancing Limited Resources under Uncertainty}
\author{Criston Hyett, Laurent Pagnier, Jean Alisse, 
Lilach Sabban,\\ Igal Goldshtein and Michael Chertkov\\
\{cmhyett,laurentpagnier,chertkov\}@math.arizona.edu,\\ 
\{Jean.Alisse,Lilah.Saban,Igal.Goldshtein\}@noga-iso.co.il}

\date{\today}

\begin{document}
\twocolumn
\thispagestyle{empty}

\maketitle

\begin{abstract}
    We build and experiment with a realistic but reduced natural gas model of Israel. The system is unusual because (a) it is controlled from a limited number of points which are at, or close to, the gas extraction sites offshore of Israel's Mediterranean coast; (b) control specifies average flux at inlet, not pressure; (c) there are no inland compressors to regulate pressure; (d) power system is the main consumer of gas (70\% of Israel's power is generated at gas-fired power plants distributed across the country). Nature of the system suggests that a special attention should be given to understanding dynamics driven by fast transients in gas consumption meeting intra-day variations in the electricity demand, and accounting for increasing role of uncertain renewable generation (mainly solar). Based on all of the above we pose and resolve a sequence of dynamic and control challenges, such as: How to time ramping up- and down- injection of gas to guarantee a healthy intra-day line-pack which meets both pressure constraints and gas-extraction patterns? We report simulation results and utilize monotonicity properties of the natural gas flows which render robustness of our conclusions to the uncertainties of the edge withdrawals of gas. 
\end{abstract}

\section{Introduction} \label{sec:intro}
This is the first manuscript of the joint NOGA (Power System Operator of Israel) and UArizona team. We aim in this manuscript
\begin{itemize}
    \item To explain specifics of the Natural Gas system of Israel, describe its minimal model and formulate main operational challenges related to uncertainty in production and consumption, also amplified by limited availability of resources. [Section \ref{sec:motivation}] 
    \item To describe basic modeling tools which allow us to study the system, including description of equations and of the software utilized. [Section \ref{sec:model}]
    \item To formulate operational scenarios, reflecting the aforementioned uncertainty, and present analysis and control solutions. [Section \ref{sec:scenarios}]
    \item To sketch further plans to extend the demonstrated methodologies towards (a) joint, and thus more realistic, modeling and control of natural gas and power systems of Israel; but also (b) formulating further practical and academic challenges towards exporting the developed methodology to other energy systems of the size comparable to the one of Israel. [Section \ref{sec:conclusion}.]
\end{itemize}

\section{Motivation, Data \& Sources} \label{sec:motivation}


Following the signing of the global agreement at the Paris Climate Conference in 2015, Israel set long-term goals to reduce greenhouse gas emissions in an effort to take part in the global action against climate change. Since the discovery of substantial offshore gas fields off the coast of Israel, Natural Gas (NG) has become the main fuel for electricity production in the country. In order to reduce Israel's carbon footprint, a decision has been made to close both major coal-fueled power plants at Hadera and Ashkelon in the future, and to convert them to gas fueled units. Apart from  renewable energy, this will render NG as the sole source of energy for electricity in the country. The Ministry of Energy plan is aimed at fulfilling Israel’s role in the agreement and at promoting an efficient, green economy. The objective is to generate up to 30\% of electricity using renewable energy by 2030, and potentially generate up to 100\% of electricity this way by 2050. The increasing share of renewable energy presents several challenges in the interim however, particularly with regard to the effects on the natural gas system, which must balance the intermittent and variable electricity production of uncontrollable renewable sources.

As of 2020, more than 50\% of Israel's electricity is produced from NG. The yearly demand for NG recently increased beyond 11 billion cubic meters (BCM). During the transition to the 2030 goal, that share is expected to increase up to 80\%-85\% in certain years. Moreover, two agreements were signed with Egypt and Jordan stipulating that Israel will export 130 BCM of NG from the Leviathan and Tamar gas fields over the next ten years. (At the moment, the non-electricity end use of NG in Israel is small, amounting to around 9000 MMBTU/h or around 10\% of the hourly typical flow).
Those data show that there are many issues surrounding the NG network in Israel over the next decade.

\href{https://www.noga-iso.co.il/en/}{NOGA}, Israel Independent System Operator (IISO), is the newly founded Israeli Electric System Operator. Its mandate is to act to ensure continuous electricity supply, at the required reliability and quality level, to all electricity consumers, in normal and emergency system conditions, and to manage the wholesale electricity market operations competitively and equitably. In addition to managing day-to-date operations, IISO is also in charge for planning development of the generation system, including recommendations for the required generation and storage capacity, maintaining optimal mix, location and timing for integration of generation and storage facilities. IISO is tasked to set up 
criteria for planning the development of the generation system. The company mandate also includes many other aspects of the transmission system planning, such as related to data forecast, transformer placement and characterization and, overall, formulation of a multi-year plan for the transmission system development. It is recognized, that to achieve all the goals IISO needs state-of-the-art tools to model the Natural Gas system and its interaction with the Electricity network.


\begin{figure}
\includegraphics[width=0.33\linewidth, valign=t]{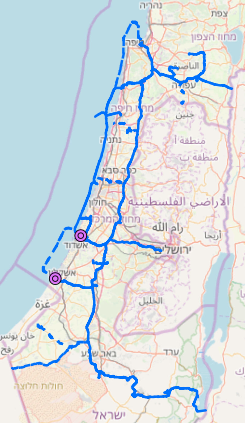}%
\includegraphics[width=0.33\linewidth, valign=t]{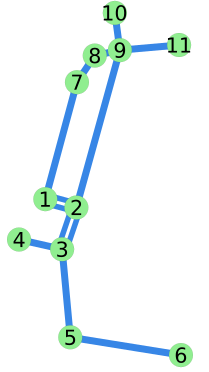}%
\includegraphics[width=0.33\linewidth, valign=t]{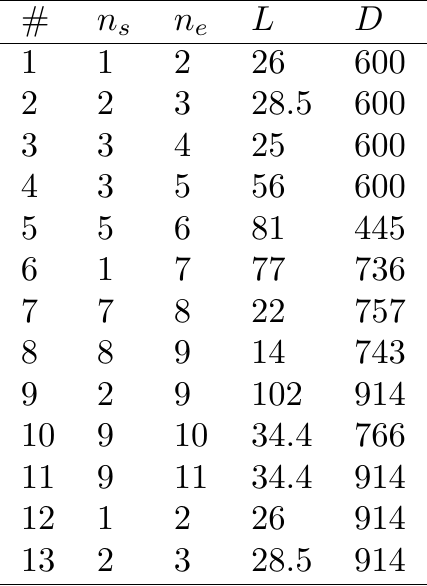}
\caption{
Description of the Israeli gas system: (left) Sketch of the true system. (center) Map of the reduced 11-node system that we use for this study. (right) List of pipelines in the simplified system. Their start node ($n_s$), end node ($n_e$), length $L$ in km, and diameter $D$ in mm. 
\label{fig:map}
}
\end{figure}

\begin{figure}
    \centering
    \includegraphics[width=0.75\linewidth]{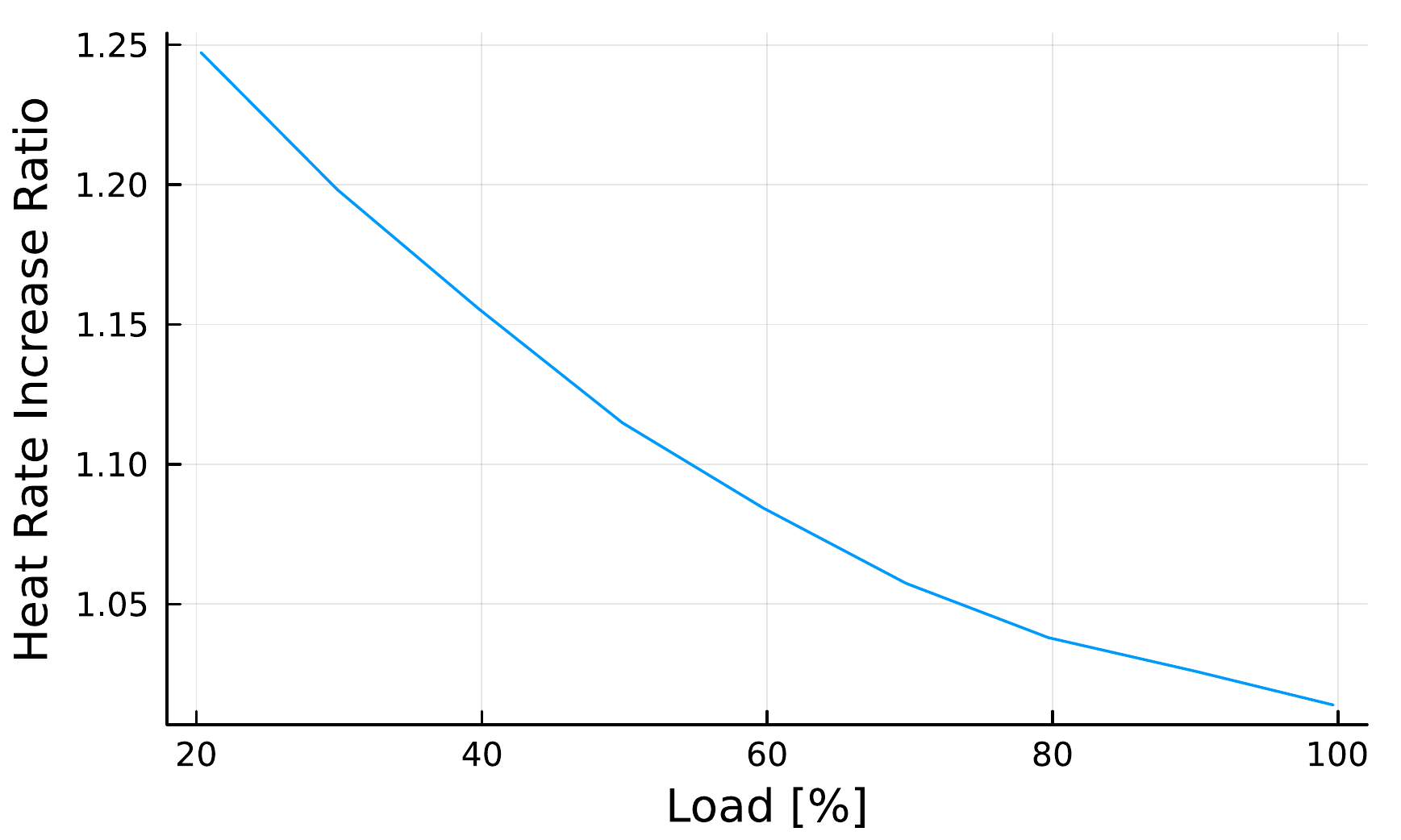}%
    \caption{Typical Gas Turbine Efficiency Curve}
    \label{fig:efficiency}
\end{figure}
In this manuscript we present a reduced, schematic gas model of Israel based on open-source public resources, shown in Figure \ref{fig:map}. The data source is the website of the Israel Natural Gas Authority www.ingl.co.il.
In order to perform the simulations, we need to create time series of the gas withdrawals at the model nodes. To do this, we use public data published by the Electricity Authority from the website \url{https://www.gov.il/he/departments/general/hovatdivuahnetunim}. The data are in the form of half-hourly time series of electricity production (in MW) at all units. We convert the power data to gas consumption with the help of typical gas turbines efficiency curves (see Figure \ref{fig:efficiency} for an example of such a curve). 

\section{Dynamic Modeling of Pipe Flow}\label{sec:model}
We utilize dynamic modeling discussed in \cite{osiadacz_simulation_1984,steinbach_pde_2007,Chertkov2015Cascading} (see also references therein)  and specific algorithmic scheme from \cite{Gyrya2019AnExplicit} called \href{https://github.com/kaarthiksundar/GasTranSim.jl}{GasTranSim.jl} implemented by Los Alamos National Laboratory team in Julia \cite{GasTranSim}. 

Basic equations, expressing conservation of mass and momentum and  describing flow of gas in a single pipe (with gas velocity much smaller than the speed of sound), are  \cite{osiadacz_simulation_1984,steinbach_pde_2007}:
\begin{align}\label{eq:mass}
    & \partial_t \rho + \partial_x \phi = 0,\\ \label{eq:momentum}
    & \partial_t \phi + \partial_x p = -\frac{\lambda}{2D} \frac{\phi |\phi|}{\rho},
\end{align}
where $\rho(t;x)$, $\phi(t;x)$ and $p(t;x)$ are the gas density, mass flux and pressure measured at the moment of time $t$ at the position $x$ along the pipe; $\lambda$ is the Darcy-Weisbach friction factor of the pipe (per diameter, $D$).
Eqs.~(\ref{eq:mass}) -(\ref{eq:momentum}) are supplemented by the equation of state, relating pressure and density
\begin{equation}
    p = Z(\rho,T)RT\rho
\end{equation}
In this study we use the CGNA formula for $Z(\rho, T)$, see e.g. \cite{menon}.
Effects of gravity and temperature variations along the system are ignored, as the effects are small and our representation is schematic
\footnote{Taking into account the temperature variations is only important in countries with sub-zero winter temperatures. 
In Israel, this is not the case and we make the assumption of a constant temperature equal to 15 Celsius degrees. Accounting for the effects of elevation/gravity is of a concern only at the node \#6, located 
in the Dead Sea area, roughly 400 meters below the sea level, and at the node \# 5, located close to  Beer-Sheva at 300 meters above the sea level. Given that ignorning this effect results only in a relatively minor pressure drop of $3$ to $5$ bars in the southern part of the system, we ignore it for now in the simplest version of the model.}.

The single pipe description extends to a system of pipes. Each pipe is characterized by three parameters: diameter, length, and the friction factor per diameter. Each node prescribes a boundary condition for one side of (at-least) one pipe at all instances in time.
Additionally, nodes and pipes are joined via condition of mass conservation (Kirchoff's rule, that the mass entering a junction must equal the mass exiting the junction).
As is standard, demand nodes specify a flux withdrawal as a function of time. The system in question however, additionally specifies flux at supply nodes. Thus all supplied boundary conditions are on flux.

Denoting $\rho_{ij}, \phi_{ij}$ to be the dynamic variables on the pipe from node $i$ to node $j$, and $\rho_n,\phi_n$ to be the density and flux at node $n$, we write the full system to be solved as:
\begin{align}
    \partial_t \rho_{ij} + \partial_x \phi_{ij} &= 0\\
    \partial_t \phi_{ij} + \partial_x p_{ij} &= -\frac{\lambda_{ij}}{2D}\frac{\phi_{ij}|\phi_{ij}|}{\rho_{ij}}\\
    \text{subject to initial}& \text{ and boundary conditions:} \nonumber\\
    \rho_{ij}(x,0) &= \rho_{0,ij}(x)\\
    \phi_{ij}(x,0) &= \phi_{0,ij}(x)\\
    \phi_n(t) &= d_n(t)\\
    \sum_{j \in \mathcal{E}} \phi_j S_{ij} + d_{j} &= 0
\end{align}

Where $S_{ij}$ is the cross-section of the pipe. Initial conditions for density and mass-flux in the system are constructed based on actual operational data. To solve for dynamics of mass flows and pressures across the system we use the staggered-grid approach of \cite{Gyrya2019AnExplicit} which is an explicit, conservative, second order, finite difference scheme, stable given a CFL condition is satisfied.
We remind that, as of now, the  Israel system does not contain compressors.

\section{Operational Scenarios} \label{sec:scenarios}

We search for robustness in the face of challenging scenarios. To this end we construct scenarios that represent two basic phenomena:
\begin{itemize}
    \item Moderate uncertainty at demand nodes, represented through addition of a random noise at the consumption site, e.g. associated with  response of gas-generators to renewable fluctuations on the electric-side of the system.
    \begin{equation}
        d_i(t) \to X_i(t)
    \end{equation}
    where
    \begin{equation}
        dX_i(t) = \alpha(d_i(t) - X_i(t)) + \gamma dW
    \end{equation}
    is a Ornstein–Uhlenbeck process - designed so that the mean is our nominal demand, $\mathbb{E}[X_i(t)] = d_i(t)$, and the variance approaches a constant exponentially fast:
    \begin{equation}
        \text{Var}(X_i(t)) = \frac{\gamma}{2\alpha}\left(1 - e^{2\alpha t} \right)
    \end{equation}
    The parameters were tuned heuristically to ensure $\alpha$ the mean was respected, and the variance approaches
    \begin{equation}
        \text{Var}(X_i(t)) \approx 0.01 \mu_i^2
    \end{equation}
    With $\mu_i$ being the mean withdrawal of node $i$. The noise for each demand is uncorrelated, a conservative approach which ignores geography and climate scales. Implementation on an actual system should be accompanied with data analysis and forecasts to determine realistic noise types for demands.
    \item Abrupt changes -- which we coin \emph{insults} -- that occur due to malfunction, weather or other exogenous circumstance. We focus particularly on supply challenges. That is, given a supply profile $s(t)$
    \begin{equation}
        s(t) \to s(t) + \Theta(t-T)\Gamma(t)
    \end{equation}
    where $\Theta$ is the Heaviside function, $T$ is the time of insult, and $\Gamma$ is the perturbation. For example, $\Gamma(t) = -s(t)$, simulates a complete loss of supply at time $T$.
\end{itemize}

In addition to studying bare ``do nothing'' scenarios we will also analyze mitigation by controls.  We assume that controls at the injection points (off-shore extraction sites at node \#1 and node \#8 in Fig.~\ref{fig:map}) and consumption sites, are step-wise. Operationally - due to the close coupling of Israel's gas extraction and delivery - step-wise control is preferable on the supply-side, while significantly idealized for demand nodes.

Focusing on the most challenging cases on fast ramps up (and, possibly, ramps down too) of consumption, e.g. around the time of sundown (or sunrise) in winter, we limit our analysis in this manuscript to \emph{prescribed} control. This allows us to conduct intuitive and easy to interpret tests of possible options. The goal of  this exploration is to analyze how the system operator can manage gas transients in line pack, meet power demand evolving throughout a day (or a number of days)  while also not exceeding the gas system pressure constraints. 

\begin{figure}
    \centering
    \includegraphics[width=0.4\textwidth]{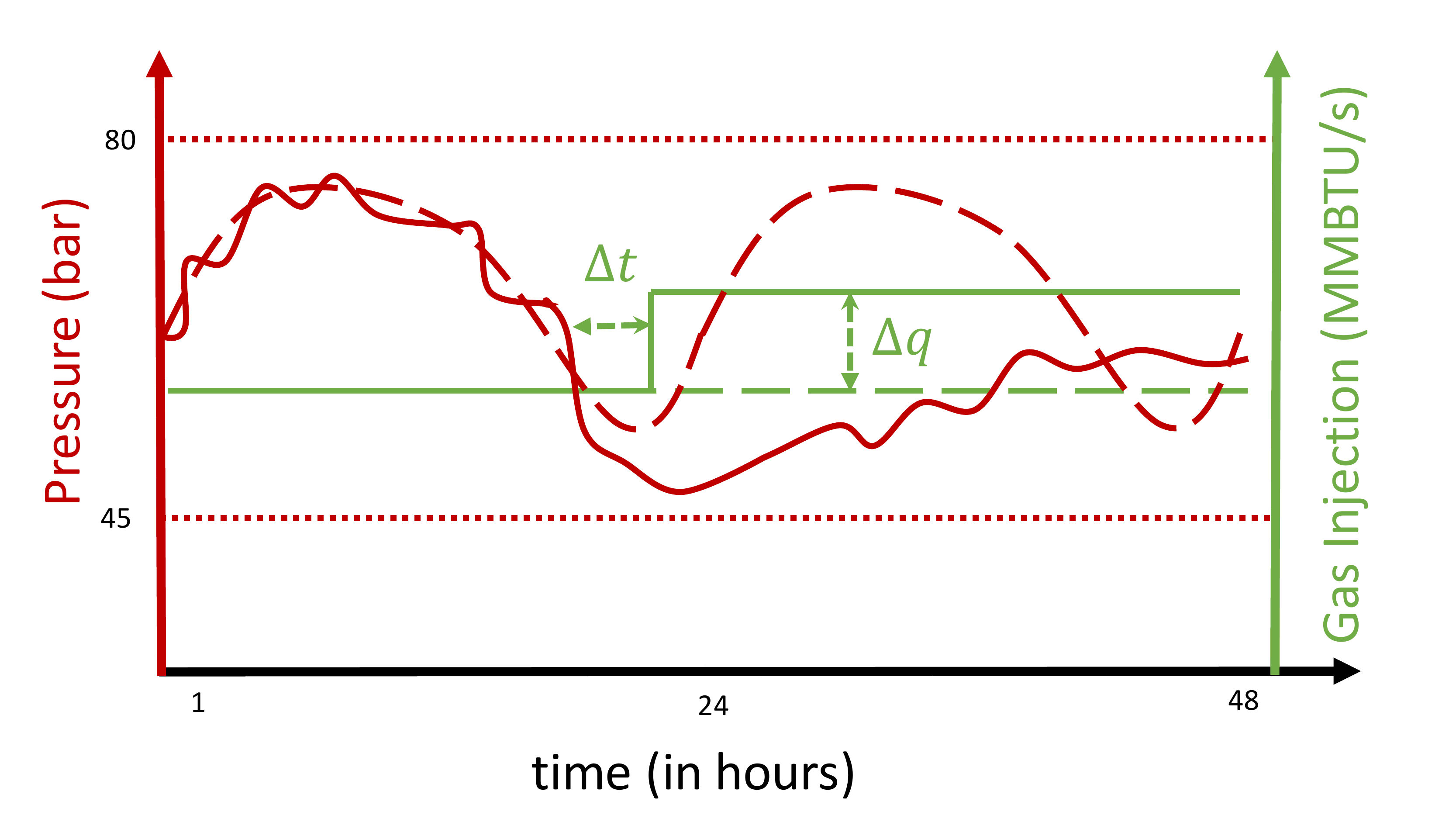}
    \caption{ Schematic illustration of the use cases of the system operations during an exemplary period of 48 hours. We show forecasts (long-dashed curves) and actual profiles (solid curves) for a pressure at a node and an injection at an entry point to the system,  where a control is applied in response to an insult. See text for details.
    \label{fig:schematic}}
\end{figure}

Our approach is illustrated schematically in Fig.~(\ref{fig:schematic}). Long-dashed green curve shows injection profile on one of the entry points to the system.  Typically,  it is flat where the value (in $MMBTU/s$) is computed based on the balanced forecast of consumption over the entire system and on the operational split of responsibilities between the injection points. Long-dashed red curve illustrates  forecast for a pressure at a node of the system.  Solid red curve shows an actual pressure profile,  which,  subject to typical uncertainty, largely follows the forecast till $\approx$ (hour) 23:00, when a significant insult occur and the deviation from the forecast becomes significant. Short-dashed red curves mark node-specific upper and lower limits for the pressure profile. Solid green curve shows operation profile at the injection point, which remains operational (in this use case) through out the 48 hours of observation, however the actual injection profile is not flat -- it follows forecast till hour $\approx$ 23:00+$\Delta t$ when a step-wise control action, responding to the $\approx$ 23:00 emergency, is applied. Selection of the proper time delay $\Delta t$, of the control response, and of the respective amplitude, $\Delta q$, constitutes a major operational challenge. 


The approach to monitoring and prescribed control of the system is summarized in Table \ref{table:scenarios}. We present six scenarios of progressive stress. Each scenario is illustrated with a figure summarizing respective dynamics. (Animations of the pressure timeseries across the network, as well as software and data to reproduce all the results included here, can be found at \url{https://github.com/cmhyett/FluxControlLinepack}.)

We define a pressure crossing as when the nodal pressure falls below 50bar. The survival time $\tau$ is defined as the time between initiation of the insult and the time to first pressure crossing \textit{at any node}. 
Note that due to integrated random fluctuations of demand, we obtain distributions of survival times, shown for example in figure[\ref{fig:scen5}] as a shaded region about the median, annotated above with the node at which the pressure crossing occurred. To keep the picture clear, the figures only show pressure crossings for a subset of our nodes, namely 9,1, and 6, as they yield information regarding the north, middle and south of our network respectively. 

\begin{enumerate}
    \item Fig.~(\ref{fig:scen1}) shows pressures (dashed) and linepack (solid), in a flux-controlled, nominal week in August. It serves as our base case, and importantly, because of the constant flux at supply nodes, the temporal variation of pressure across the network is governed by the intra-day demand fluctuations.
\begin{figure}
    \centering
    \includegraphics[width=0.47\textwidth]{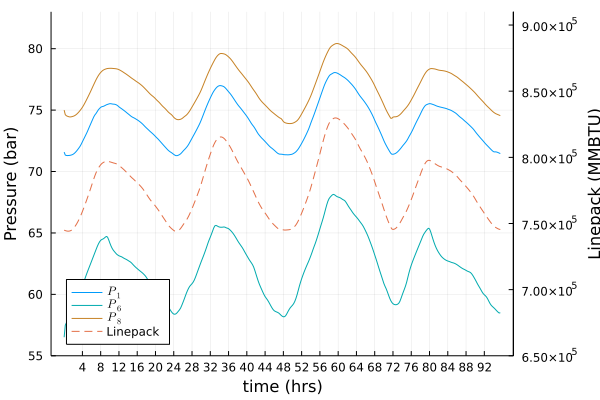}
    \caption{Nominal week in August, with no uncertainty.}
    \label{fig:scen1}
\end{figure}

    \item Fig.~(\ref{fig:scen2}) adds random fluctuations to demands on this nominal week - modeling uncertainty from exact power demand and generation of renewables. For each node, we add noise distributed uniformly with width of 5\% of nominal demand at that node. These small perturbations integrate over time, leading to significant linepack and pressure deviations from the mean.
    \begin{figure}
    \centering
    \includegraphics[width=0.47\textwidth]{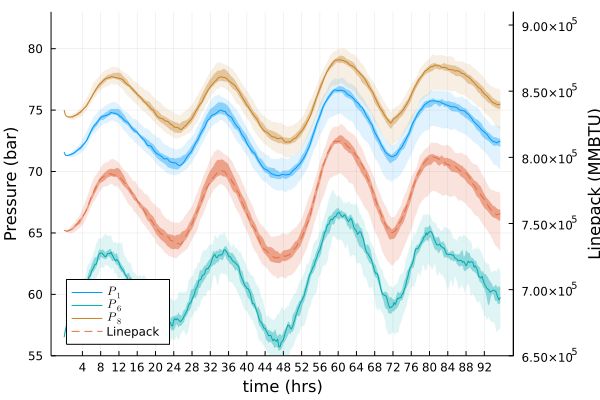}
    \caption{Nominal week in August, with empirical noise added to demand curves. Notice the drift of linepack and pressure jitter, consistent with what was predicted in \cite{chertkov_cascading_2014}. Using a Monte-Carlo with 50 simulations, we plot filled regions of containing the middle 75\%, 25\%, and median using increasing color intensities.}
    \label{fig:scen2}
\end{figure}

    \item Fig.~(\ref{fig:scen3}) introduces an ``insult'' indicating off-nominal or emergency operation.  Particularly, supply at node \#1 (one of our two supply nodes) are closed at hour 36 in the simulation. We continue to run the simulation to observe the rate of linepack decay and the sequence of pressure crossings.
    The survival time in this scenario is
    \begin{equation}
        \tau = 4.13 \pm 0.38 \text{ hrs}
    \end{equation}
    where the $0.38$ is the standard deviation.
    This information can be translated to spatiotemporal ``vulnerability'' of the gas network.
\begin{figure}
    \centering
    \includegraphics[width=0.47\textwidth]{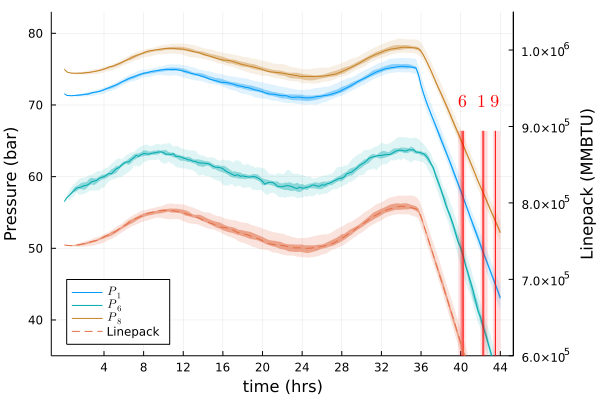}
    \caption{Scenario 3 results, an insult at a crest of the linepack at $t=36\text{hrs}$. Using a Monte-Carlo with 50 simulations, we plot filled regions containing the middle 75\%, 25\%, and median using increasing color intensities.}
    \label{fig:scen3}
\end{figure}

    \item Fig.~(\ref{fig:scen4}) introduces the same insult as scenario 3, but at hour 48 instead of hour 36. In particular, this corresponds to the insult occurring at a trough of the linepack curve instead of a peak. We highlight first that the time to first pressure crossing survival time is shorter in this case
    \begin{equation}
        \tau = 3.58 \pm 0.89 \text{ hrs}
    \end{equation}
    This simple statement that the survival time of the network depends on the start time of an insult is the result of complicated interactions between demand-node boundary conditions as well as network topology.
\begin{figure}
    \centering
    \includegraphics[width=0.47\textwidth]{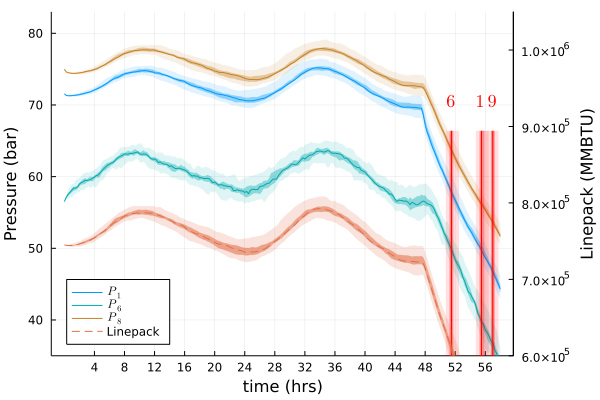}
    \caption{Scenario 4 results, an insult at a trough of the linepack at $t=48\text{hrs}$. Using a Monte-Carlo with 50 simulations, we plot filled regions containing the middle 75\%, 25\%, and median using increasing color intensities.}
    \label{fig:scen4}
\end{figure}

    \item Fig.~(\ref{fig:scen5}) begins introducing control, attempting to mimic "human in the loop" control of the network under the insult described in scenario 4. In this scenario, the operator at the remaining supply (node 8), increases to the max flow-rate a half hour after the insult begins ($t=48.5\text{hrs}$). This control stabilizes the linepack in the short-term, but fails to handle the daily ramp near hour 60. In particular, the south of the network, far from the remaining supply contains all of the pressure crossings.
    \begin{equation}
        \tau = 14.17 \pm 4.07 \text{ hrs}
    \end{equation}
\begin{figure}
    \centering
    \includegraphics[width=0.47\textwidth]{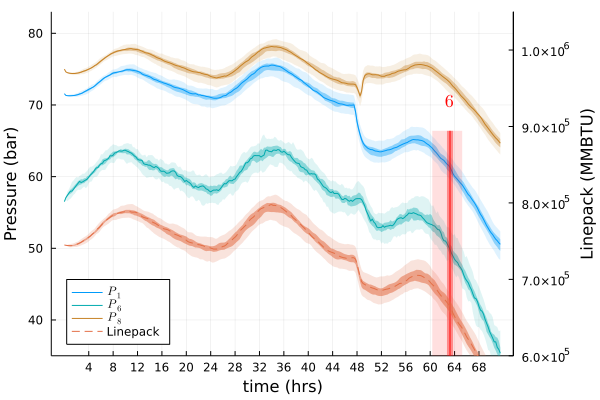}
    \caption{Scenario 5 results, introduces a step-wise increase in supply at node 8 half an hour after the insult ($t=48.5\text{hrs}$).}
    \label{fig:scen5}
\end{figure}

    \item Fig.~(\ref{fig:scen6}) finally builds upon scenario 5 to additionally curtail demand 2 hours after the insult ($t=50\text{hrs}$). This translates to a variety of potential action, from high penalty demand-response (as demonstrated during heat waves in California) to the transition of natural gas plants to alternative fuels, or the utilization of stored power. 
    \begin{figure}
    \centering
    \includegraphics[width=0.47\textwidth]{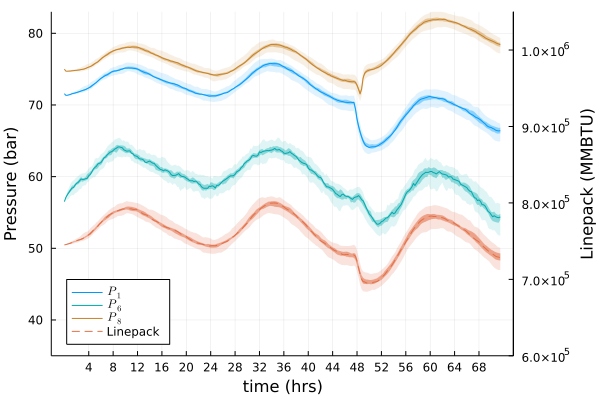}
    \caption{Scenario 6 results, curtails demand at $t=50\text{hrs}$, to maintain minimum pressures across the network.}
    \label{fig:scen6}
\end{figure}

\end{enumerate}

\section{Monotonicity}\label{sec:monotonicty}
In our scenario selection, we were intentionally coarse, preferring severe and abrupt challenges, while searching for the mildest controls for remedies. This is intentional, as previous work has given us monotonicity guarantees \cite{Vuffray2015Monotonicity,zlotnik_monotonicity_2016}. That is, we have that for any less-severe challenge (as in the more physically realistic scenario of a slow ramp down of a supply instead of our simulated near-immedate shut off), our pressures will be bounded below by the most severe case, and thus in turn our estimates for survival times are in fact lower, conservative bounds.

We advocate this approach for schematic expositions as it is only with full operational and procedural knowledge that one can obtain tight estimates - that is, work inherently reliant on proprietary data.

It should also be noted that monotonicity can yield bounds for the more usual scenario of pressure and linepack drift resulting from integration of stochasticity due to renewables, such drift can be seen in  Fig.\ref{fig:scen2}. However, monotonicity bounds were derived without relation to probability, thus it is likely that future work can tighten these bounds, avoiding expensive simulation except when full distributional knowledge is needed.

\section{Conclusion}\label{sec:conclusion}
We investigate the spatiotemporal response of a reduced model of Israel's NG network to prescribed insults and human-in-the-loop controls in order to evaluate robustness and suggest control strategies. To reiterate, Israel's network is unique because of the absence of a compressor, and that the inlets specify flux, not pressure. Further, we perform this study looking towards the increased importance of NG to mitigate increasing stochasticity in demands expected in the coming years as coal is phased out, and renewables grow.

The specification of flux vs pressure leads to the pressure timeseries of the network being dominated by daily demand curves as shown in Fig.~(\ref{fig:scen1}), increasingly susceptible to pressure drift from stochastic fluctuations in nominal demands.

Further, we call out the importance of robustness of the network not simply to insults, but to insults at any time - leading to the idea of "system reserve" being time and spatially dependent.

Future work will improve on modeling to more completely capture uncertainty propagation through the network, and its influence and interaction with control strategies. We envision extending the prescribed control, also reinforced by monotonicity \cite{Vuffray2015Monotonicity,zlotnik_monotonicity_2016}, developed in this manuscript with the powerful optimization approaches developed to account for dynamic optimization over compressors \cite{rachford_optimizing_2000,carter_optimizing_2003,rachford_using_2009,zlotnik_optimal_2015,zlotnik_using_2016}, e.g. to evaluate benefits of adding compressors to the NG system of Israel. We also plan to carry on a comprehensive modeling and control of the combined power and gas system of Israel, in the spirit of the approach highlighted in \cite{carter_impact_2016,Zlotnik2017Coordinated}. 

\section*{Acknowledgments}
The authors are grateful for discussions with Anatoly Zlotnik, Vitaliy Gyrya, and Kaarthik Sundar.
CH, LP and MC acknowledge support from UArizona, as well as travel grant from NOGA Israel.
CH acknowledges support from NSF's Data-Driven Research Training Grant.

\bibliographystyle{ieeetr}
\bibliography{gas.bib,MishaPapers.bib,gas1.bib}

\section*{Appendices}\label{sec:app}
\begin{table*}
    \centering
    \begin{tabular}{|p{0.1\textwidth} | p{0.4\textwidth} | p{0.4\textwidth}|}
    \hline
    Scenario \#  &   \multicolumn{1}{c|}{Description} &   \multicolumn{1}{c|}{Features}\\
    \hline
    \multicolumn{1}{|c|}{1} & A reference week in August. & Pressure variation in flow-control regime \\ 
    \hline
    \multicolumn{1}{|c|}{2} & Scenario \#1 with empirical noise added to demand curves, supplies unchanged. & Linepack and pressure drift when using flow-control with uncertain demand.\\ 
    \hline
    \multicolumn{1}{|c|}{3} & Scenario \#2 with insult at node 1. & Introduce the notion of survival time, and set baseline without any controls.\\ 
    \hline
    \multicolumn{1}{|c|}{4} & Scenario \#3 with insult time change to trough of linepack timeseries. & Illustrate that survival times change with timing of insult.\\ 
    \hline
    \multicolumn{1}{|c|}{5} & Scenario \#4 with step-wise supply increase from node \# 8. & Survival times lengthen, but become less certain. \\
    \hline
    \multicolumn{1}{|c|}{6} & Scenario \#5 with step-wise curtailing of demand. & No low pressure crossings are found. The high pressure at node \# 8 shows need for finer control.\\
    \hline
    \end{tabular}
    \caption{Description of scenarios.}
    \label{table:scenarios}
\end{table*}

\section*{Author Biographies}

\textbf{Criston Hyett} - is a Ph.D. student in Applied Mathematics at the University of Arizona. He is advised by Misha Chertkov, and interested in data-enhanced dynamical systems modeling, encompassing physics-informed machine learning, reduced order modeling and uncertainty quantification.

\noindent
\textbf{Laurent Pagnier} - is a Visiting Assistant Professor at the University of Arizona. His main research interest is the application of Machine Learning techniques to the operation of large infrastructures. He is particularly interested in reinforcing the interpretability and trustworthiness of ML methods which are paramount to increase their acceptance and usage by practitioners.

\noindent
\textbf{Jean Alisse} - received his Phd in Physics in 1999 (Paris University). He worked during 11 years at the Israel Electric Company in the Planning and Division company.
There he headed a Modeling group which focused on CFD problems and Natural Gas dynamics. Since the creation of Noga, the Israel Independent System Operator, he has been pursuing the same works, with a stress on developing models for Gas dynamics simulation and Gas-Power systems. 

 \noindent
\textbf{Lilach Sabban} - received her PhD from the Technion Israel Institute of Technology. Lilach specializes in fluid dynamics as a mechanical engineer. She is a researcher at Noga, Israel Independent System Operator,  involved in a range of renewable energy projects and has been investigating natural gas dynamics.

\noindent
\textbf{Igal Goldshtein} - received his B.Sc from the Technion Israel Institute of Technology. He worked at Israel electric company for 12 years, seven of them as Gas-Turbine dispatcher in the system operator control center. During the last two years, at Noga, the Israel System Independent Operator, he handles the Operator Training System (OTS). He also researches operation of Israel electric system under various stress conditions.

\noindent
\textbf{Michael (Misha) Chertkov} - is Professor of Mathematics and chair of the Graduate Interdisciplinary Program in Applied Mathematics at the University of Arizona since 2019. He focuses in his research on foundational problems in mathematics and statistics applied to physical systems, in particular fluid mechanics, to engineered systems such as energy grids, and to some bio-social systems. Dr. Chertkov received his Ph.D. in physics from the Weizmann Institute of Science in 1996, spent three years at Princeton University as a R.H. Dicke Fellow in the Department of Physics, and joined Los Alamos National Laboratory in 1999, initially as a J.R. Oppenheimer Fellow and then as a Technical Staff Member in Theory Division. He has published more than 250 papers, is a fellow of the AAAS, a fellow of the American Physical Society and a senior member of IEEE.

\end{document}